\documentclass[a4paper, 12pt]{article}
\usepackage{makeidx}
\usepackage{enumerate, theorem}
\usepackage{amsmath, amsfonts, amssymb}
\usepackage[height=22.5cm, width=15cm]{geometry}
\usepackage[all]{xy}
\usepackage{mathrsfs, yfonts}
\usepackage{graphicx}
\DeclareMathOperator{\C}{\mathbb{C}}
\newcommand{\parag}[1]{\paragraph{\sc{#1.}} }

\newtheorem{thm}{Theorem}[subsection]

\newtheorem{cor}[thm]{Corollary}
\newtheorem{prop}[thm]{Proposition}
\newtheorem{lemma}[thm]{Lemma}

\begin{document}

\title{Note on Lisbon integrals and their associated $D-$modules}

 \author{Daniel Barlet\footnote{Institut Elie Cartan, G\'eom\`{e}trie,\newline
Universit\'e de Lorraine, CNRS UMR 7502   and  Institut Universitaire de France.}.}

\maketitle

\parag{Abstract} The aim of this Note is to specify the links between the three kinds of Lisbon integrals, trace functions and trace forms with the corresponding $D-$modules.\\

\parag{AMS Classification} 44A99, 32C35, 35A22, 35 A27, 58J15.

\tableofcontents

\section{Introduction} 
The aim of this Note is to specify the links between the three kinds of Lisbon integrals, trace functions and trace forms and the corresponding $D-$modules defined by the systems $S\, 1), S\, 2)$ and $S\, 3)$ (see below). The new results here are given by  theorem \ref{solutions} (in fact in theorem \ref{caract.}) which gives a characterization of trace forms as solutions of the system $S\, 2)$ and its two corollaries \ref{charact. bis} and \ref{3 to 1}. The first  corollary gives a simpler characterization of trace functions and the second specifies the link between the system $S\, 1)$ and $S\, 3)$.
 We conclude by showing that there is a simple linear system relating the interpolation polynomial of an entire function $f \in \mathcal{O}(\C)$ with the solution of the system $S\, 3)$ given by the vector Lisbon integral associated to $f$.

\section{Lisbon integrals, Traces and System of Partial Differential Equations}

\subsection{Definitions and results}

\parag{Notation} For $\sigma := (\sigma_1, \dots, \sigma_k) \in \C^k$ we denote $P_\sigma(z) := \sum_{h=0}^k (-1)^h.\sigma_h.z^{k-h} $ with the convention $\sigma_0 \equiv 1$.\\

For basic properties on $D-$modules the reader may consult \cite{[Bjork]}.\\

We shall use the following three kinds of  Lisbon integrals
\begin{enumerate}[L 1)]
\item The first kind, which  associates to $f \in \mathcal{O}(\C)$ the holomorphic function on $N := \C^k$ given by:
$$ F(\sigma) := \frac{1}{2i\pi}\int_{\vert \zeta\vert = R} f(\zeta).\frac{P'_{\sigma}(\zeta)}{P_{\sigma}(\zeta)}.d\zeta \quad {\rm for} \   R \gg \vert\vert \sigma \vert\vert.$$
\item The second kind which associates to $f \in \mathcal{O}(\C)$ the holomorphic function on $N := \C^k$ given by
$$ \tilde{F}(\sigma) := \frac{1}{2i\pi}\int_{\vert \zeta\vert = R} f(\zeta).\frac{d\zeta}{P_{\sigma}(\zeta)} \quad {\rm for} \   R \gg \vert\vert \sigma \vert\vert.$$
\item The vector kind which associates to $f \in \mathcal{O}(\C)$ the holomorphic function $$\Phi := \begin{pmatrix} \varphi_0 \\ \varphi_1\\ . \\ . \\ \varphi_{k-1} \end{pmatrix} $$  on $N := \C^k$ with values in $\C^k$ given by
$$ \varphi_h(\sigma) := \frac{1}{2i\pi}\int_{\vert \zeta\vert = R} f(\zeta).\frac{\zeta^h.d\zeta}{P_{\sigma}(\zeta)} \quad {\rm for} \   R \gg \vert\vert \sigma \vert\vert$$
where $h$ is in $[0, k-1]$.
\end{enumerate}
We shall also consider the following three cases of trace:
\begin{enumerate}[T 1)]
\item The trace function  associated to $f \in \mathcal{O}(\C)$ which  is the holomorphic function on $N := \C^k$ given by:
$$ T(f)(\sigma) = \sum_{P_\sigma(z_j) = 0} f(z_j) .$$
\item The trace form\footnote{We leave out the holomorphic volume form $d\sigma$ on $N$, the holomorphic volume form  $dz$ on $\C$ and also the $\C-$ relative  volume form $d\sigma_1\wedge \dots \wedge d\sigma_{k-1}$  of the projection $p : H \to \C$ given by $p(\sigma, z) = z$ where $H$ is the hypersuface $H := \{(\sigma, z) \in N \times \C \ / \  P_\sigma(z) = 0  \}.$}  associated to $f \in \mathcal{O}(\C)$ which is the holomorphic function on $N := \C^k$ given by:
$$ \tilde{T}(f)(\sigma) = \sum_{P_\sigma(z_j) = 0} \frac{f(z_j)}{P'_\sigma(z_j)} .$$
\item The vector trace form which is associated   to $f \in \mathcal{O}(\C)$ which is the holomorphic function $V\tilde{T}(f)$  on $N := \C^k$ with values in $\C^k$ given by
$$ V\tilde{T}(f)_h(\sigma) := \sum_{P_\sigma(z_j) = 0} \frac{z_j^h.f(z_j)}{P'_\sigma(z_j)} \quad {\rm for} \quad h \in [0, k-1].$$
\end{enumerate}
Remark that $V\tilde{T}(f)_h = \tilde{T}(z^h.f) $ for each $h \in [0, k-1]$.\\

We shall also consider three systems of partial differential equations on $N := \C^k$:
\begin{enumerate}[S 1)]
\item The first one is given by the left ideal  $\mathcal{I}$ in $D_N$ generated by the following partial differential operators :
$A_{p, q} := \partial_{p}\partial_{q} - \partial_{p+1}\partial_{q-1}$ for $p \in [1, k-1]$ and $q \in [2, k]$ and by 
$T^{m} := \partial_{1}\partial_{m-1} + \partial_{m}E \quad {\rm for}  \ m \in [2, k] $, where $E := \sum_{h=1}^k \sigma_h.\partial_h$.
\item The second one is given by the left ideal  $\tilde{\mathcal{I}}$ in $D_N$ generated by the following partial differential operators :
$A_{p, q} $  for $p \in [1, k-1]$ and $q \in [2, k]$ and by  $\tilde{T}^{m} := T^{m} + \partial_{m}$ \ for \  $m \in [2, k] $.
\item The last one is given by the following vectorial equations, where $\Phi \in (\mathcal{O}_N)^k$ is the unknown:\\
\begin{equation*}
 (-1)^{k-h}\frac{\partial \Phi}{\partial \sigma_h} = \frac{\partial (A^{k-h}.\Phi)}{\partial \sigma_k}  \quad {\rm for} \ h \in [1, k-1] \tag{$S\, 3)$}
 \end{equation*}
where $A $ is the matrix given by
$$ A :=  \begin{pmatrix} 0 & 1 & 0 & \dots & \dots & 0 \\ 0 & 0 & 1 & \dots &\dots & 0 \\ \dots & \dots & \dots & \dots & \dots &\dots \\ \dots & \dots & \dots & \dots & 1 & 0  \\  0 & 0 & \dots & \dots & 0 & 1 \\ \tilde{\sigma}_k & \tilde{\sigma}_{k-1} &\dots & \dots &\dots & \tilde{\sigma_1} \end{pmatrix} $$
\end{enumerate}
and where $\tilde{\sigma}_h := (-1)^{h-1}. \sigma_h$\  for $h \in [1, k]$.\\

It is an easy exercise on the residue formula to show that $L\, i)$ is equivalent to $T\, i)$ for $i = 1, 2, 3$, corresponding to the equalities
$$ F = T(f), \quad \tilde{F} = \tilde{T}(f) \quad {\rm and} \quad \Phi = V\tilde{T}(f) $$
when $F, \tilde{F}$ and $\Phi$ are associated to $f \in \mathcal{O}(\C)$.\\
These formulas show the holomorphy of these functions on all $N$, despite the denominators in $T\, 2)$ and $T\, 3)$.\\
It is also an easy exercise to verify that for $ R \gg \vert\vert \sigma \vert\vert$ the function $\zeta \mapsto Log(P_\sigma\big(\zeta)\big/\zeta^k\big)$ is well  defined  around the circle $\{\vert \zeta\vert = R\}$  and that for any $f \in \mathcal{O}(\C)$ we have the formula
\begin{equation*}
F(\sigma) = -\frac{1}{2i\pi}\int_{\vert \zeta \vert = R} f'(\zeta).Log\big(P_\sigma(\zeta)\big/\zeta^k\big).d\zeta \ + k.f(0) \tag{M}.
\end{equation*}
when $R$ is large enough compared with $\vert\vert \sigma \vert\vert$.
This formula will be very useful to derive the trace function $F$.\\

We summarize now the main results linking Lisbon integrals, traces and the differential systems defined above.

\begin{thm}\label{solutions}
The global solution of $S\, i)$ are exactly the traces $T\, i)$ which coincide with the Lisbon integrals $L \,i)$  for $i = 1, 2, 3$.
\end{thm}

For $i = 1$ this is proved in \cite{[B.19]}, for $i = 2$ this is proved in section 3.2  (see the theorem \ref{caract.}) and for $i = 3$ this is proved in \cite{[B-MF]} where the system $S\, 3)$  is described with functorial operations on $D-$modules.

\bigskip

We shall now give some corollaries of the theorem \ref{solutions} 

\begin{cor}\label{charact. bis}
The trace functions $F$ are the only global solutions of the restricted system 
\begin{equation*} T^m(F) = 0 \quad \forall m \in [2, k] \tag{$S\, 0)$}
\end{equation*}
\end{cor}

This corollary will be proved in section 3.3.

\begin{cor}\label{3 to 1}
 If $\Phi$ is  the solution of $S\, 3)$  associated to $f$ then the vector form $dA.\Phi$ is $d-$closed and has only its last component which is non zero. This last component is the differential  of  the trace function associated to a primitive of $(-1)^k.f$. 
  \end{cor}
  
  This corollary is also proved in section 3.3.
 
\parag{Remarks}

\begin{enumerate}
\item The function  $F$ such that $dF$ is the last component of $dA.\Phi$ is defined up to a constant. So is the primitive of $f$.
\item The converse in terms of Lisbon integrals,  so the fact that for any trace function $F$ the vector function $\Phi$ defined by $\varphi_{k-h} = (-1)^{h-1}.\partial_hF$ is solution of $S\, 3)$ is an easy consequence of the formula $(D)$ below (see the proof of point 7. of \ref{direct relations}) using the theorem \ref{solutions}.
\end{enumerate}

\subsection{Some direct links} 

Before giving the proofs of the previous results in section 3, we shall give in the next proposition direct proofs of some easy links between the Lisbon integrals, the trace functions and the three $D-$modules introduced above which do not use the theorem \ref{solutions}.

\begin{prop}\label{direct relations}
\begin{enumerate}
\item Let $F$ be the Lisbon integral associated to $f \in \mathcal{O}(\C)$. Then $ (-1)^{h-1}.\partial_h(F)(\sigma) = V\tilde{T}(f')_{k-h}$ where $V\tilde{T}(f')$  is the Lisbon vector integral associated to $f'$.
\item Also, if $F$ is the trace function associated to $f$, $(-1)^{h-1}.\partial_h(F)$ is the trace form associated to the entire function $z \mapsto z^{k-h}.f'(z)$.
\item If $F$ is a solution of $S\, 1)$ and $h \in [1, k]$ then $\partial_h(F)$ is solution of $S\, 2)$.
\item If $\Phi$ is  a solution of $S\, 3)$, then $A.\Phi$ is also a solution of $S\, 3)$.
\item Also if $\Phi$ is the vector Lisbon integral associated to $f$, then $A.\Phi$ is the vector Lisbon integral associated to $z \mapsto z.f(z)$.
\item The solutions of the systems $S\, 1)$ and $S\, 2)$  are stable by the action of the vector fields 
 $U_0 := \sum_{h=1}^k h.\sigma_h.\partial_h$ and  $U_{-1} := \sum_{h= 0}^{k-1} \ (k-h).\sigma_h.\partial_{h+1} $.
\item Also  Lisbon integrals of type $S\, 1)$ and $S \, 2)$ are  stable by $U_0$ and  $U_{-1}$.
\item If $\Phi$ is solution of $S\, 3)$ then the vector form $\omega :=  dA.\Phi$ is $d-$closed and has only its last component which is non zero. This implies that there exists a holomorphic function  $F$ on $N$ such that $(-1)^{h-1}.\partial_h(F) = \varphi_{k-h} \quad \forall  h \in [1, k]$. 
\end{enumerate}
\end{prop}

Note that the fact that in point 8. above the function $F$ is a trace function is proved in corollary \ref{3 to 1}. But it is not obvious to give a direct and easy proof without using the theorem \ref{solutions}.

\parag{Proof} The points 1. and 2. are obtained by derivation of the formula $(M)$ above.\\
The point 3. is consequence of the commutations relations
$$ [A_{p, q}, \partial_h] = 0 \quad {\rm and} \quad [T^m, \partial_h] = -\partial_h\partial_m \quad {\rm which \ implies} \quad \tilde{T}^m.\partial_h = \partial_h.T^m .$$
To prove point 4. let $\Phi$ be a solution  of $S \, 3)$. Then we have
\begin{align*}
& (-1)^{k+h}.\partial_h(A.\Phi) = (-1)^{k+h}.\partial_h(A).\Phi + (-1)^{k+h}.A.\partial_h(\Phi) \\
& (-1)^{k+h}.\partial_h(A.\Phi) = (\partial_kA).A^{k-h}.\Phi + A.(\partial_k(A^{k-h}.\Phi)) 
\end{align*}
and using  the  relation\footnote{See lemma 2.5 in \cite{[B-MF]}}
 \begin{equation*}
 (-1)^{k+h}.\partial_h(A) = (\partial_kA).A^{k-h} 
 \end{equation*}
 and so we obtain
$$  (-1)^{k+h}.\partial_h(A.\Phi) = \partial_k(A^{k-h}.(A.\Phi)) .$$

For point 5. remark that we have $A.E(z) = z.E(z) + P_\sigma(z).V$ where $V := \begin{pmatrix} 0 \\ .\\ . \\ 0 \\ 1\end{pmatrix}$. The conclusion follows from the vanishing of the integral
$$ \frac{1}{2i\pi}\int_{\vert \zeta \vert = R}  f(\zeta).V.d\zeta .$$ 
The point 6. for $S\, 1)$ and $S\, 2)$  is an easy consequence of the commutation relations given in the lemma \ref{commutators} below.\\
To show the point 7.  first use the formula $(M)$ to obtain, if $F$ is the Lisbon integral of first kind associated to $f$:
\begin{equation*}
 \partial_h(F)[\sigma] = \frac{-1}{2i\pi}\int_{\vert \zeta\vert = R} f'(\zeta)\frac{(-1)^h.\zeta^{k-h}.d\zeta}{P_\sigma(\zeta)} .  \tag{D}
 \end{equation*}
This gives
$$ U_0(F)[\sigma] = \frac{-1}{2i\pi}\int_{\vert \zeta\vert = R} f'(\zeta).\big(\sum_{h=1}^k \ (-1)^h. h.\sigma_h.\zeta^{k-h}\big).\frac{d\sigma}{P_\sigma(\sigma)} .$$
Using now the identity
$$ \zeta.P'_\sigma(\zeta) = k.P_\sigma(\zeta) - \sum_{h=1}^k  \ (-1)^h.h.\sigma_h.\zeta^{k-h} $$
we obtain
$$ U_0(F)[\sigma] = \frac{-1}{2i\pi}\int_{\vert \zeta\vert = R} f'(\zeta).d\zeta + \frac{1}{2i\pi}\int_{\vert \zeta\vert = R} f'(\zeta).\frac{\zeta.P'_\sigma(\zeta)}{P_\sigma(\zeta)}.d\zeta  = T(z.f')(\sigma)$$
proving  the stability of $L\, 1)$ by $U_0$.\\
Also, using again the formula $(D)$ we obtain
$$ U_{-1}(F)[\sigma] = \frac{-1}{2i\pi}\int_{\vert \zeta\vert = R} f'(\zeta).\big( \sum_{h=0}^{k-1} (-1)^{h+1}(k-h).\sigma_h.\zeta^{k-h-1} \big).\frac{d\zeta}{P_\sigma(\sigma)} $$
and so
$$U_{-1}(F)[\sigma] = \frac{1}{2i\pi}\int_{\vert \zeta\vert = R} f'(\zeta).\frac{P'_\sigma(\zeta)}{P_\sigma(\zeta)}.d\zeta = T(f')[\sigma] $$
where $T(f')$ is the trace of the derivative of $f$.\\
The stability of $L\, 2)$ by the actions of $U_0$ and $U_{-1}$ is then consequence of  the formula $(D)$ and of the commutations relations
$$ [U_0, \partial_h] = -h\partial_h \quad \forall h \in [1, k], \quad [U_{-1}, \partial_h] = -(k-h)\partial_{h+1} \quad \forall h \in [1, k-1] \quad {\rm and} \quad [U_{-1}, \partial_k] = 0 .$$
To prove point 8. remark first that the relation $d\omega = 0$ is equivalent to the vanishing of
$$ -\partial_p\varphi_{k-h} + \partial_h\varphi_{k-p} \quad \forall h, p \in [1, k] $$
because the only non zero line of $dA$ is the last one which is equal to 
$$\big((-1)^{k-1}.d\sigma_k, (-1)^{k-2}.d\sigma_{k-2}, \dots, d\sigma_1\big).$$ 
So only  the last component of $dA.\Phi$ is non zero and it  is equal to $\sum_{h=1}^k (-1)^{k-h-1}\varphi_{k-h}.d\sigma_h $.\\
But we have
\begin{align*}
& (-1)^{k+h}.\partial_p\varphi_{k-h} = (-1)^{p+h}\big[\partial_k(A^{k-p}.\Phi)\big]_{k-h} \\
&  (-1)^{k+p}.\partial_h\varphi_{k-p} = (-1)^{p+h}\big[\partial_k(A^{k-h}.\Phi)\big]_{k-p}
\end{align*}
where we note $[W]_d$  the $d-$th component of the vector $W$.
\smallskip

Then, it is enough to prove that the $(k-h)-$th line of the matrix $A^{k-p}$ is equal to the $(k-p)-$th line of the matrix $A^{k-h}$. This is explained in the next lemma. $\hfill \blacksquare$\\

 \begin{lemma}\label{calcul matriciel}
 Let $P(z) := z^{k} + \sum_{h=1}^{k} (-1)^{h}.\sigma_{h}.z^{k-h}$ and consider in the algebra $R := \C[\sigma, z]\big/(P)$ the endomorphism $\mathcal{A}$ of multiplication by $z$. Define
 $$ \mathcal{A}^{j}[z^{a}] = z^{a+j} = \sum_{b=0}^{k-1} a_{a+j, b}.z^{b} .$$
 Then the matrix $\Gamma_{j} $  of $\mathcal{A}^{j}$ in the $\C[\sigma]-$basis $1, z, \dots, z^{k-1}$ of $R$ is given by
   $$\Gamma_{j} := \big(\gamma(j)_{u, v} \big)= \big(a_{v+j, u}\big).$$
  Then the $(k-q)-$th line of $\Gamma_{k-p}$ is equal to the $(k-p)-$th line of $\Gamma_{k-q}$ for each $p, q$ in $[1, k]$.
 \end{lemma}
 
\parag{ Proof} Both are equal to $(a_{2k-(p+q), 1}, \dots, a_{2k-(p+q), k-1})$.$\hfill \blacksquare$\\

\parag{Remark} The matrix of $\mathcal{A}$ in the  $\C[\sigma]-$basis $1, z, \dots, z^{k-1}$ of $R$  is the matrix $A$ introduced above.

\section{The characterization of trace forms}

\subsection{Preliminaries}

We have given an explicit system (which is $S\, 1)$) of global partial differential equations on $N = \C^{k}$ characterizing the trace functions in \cite{[B.19]}. We have also obtained an analogous system (which is $S\, 2)$) such that the trace forms are solutions of it. But it was not clear that we obtain a characterization of trace forms with it. The aim of this first section is to prove that this is also true. As a by product we obtain a characterization of the trace functions with the much smaller system $S\, 0)$.

\parag{Notations} Let $N := \C^{k}$ and $D_{N}$ be the sheaf of holomorphic differential operators on $N$. We shall also consider the Weyl algebra $W_{2} := \C[\sigma_{1}, \dots, \sigma_{k}]\langle \partial_{1}, \dots, \partial_{k} \rangle$ of global algebraic sections of $D_{N}$ where $\sigma_{1}, \dots, \sigma_{k}$ are the coordinates on $N$.\\
We shall consider the left ideals $\mathcal{I}$ and $\tilde{\mathcal{I}}$ in $D_{N}$ generated respectively by the systems $S \,1)$ and $S\, 2)$.\\

Let $\alpha$ and $\beta$ be in $\mathbb{N}^{k}$. For $\sigma^{\alpha}.\partial^{\beta}$ in $W_{2}$ we  define the weight of this element as $w(\alpha) - w(\beta) \in \mathbb{Z}$, where we put $w(\gamma) := \sum_{h=1}^{k} \  h.\gamma_{h}$ for $\gamma \in \mathbb{N}^{k}$.\\
We say that $P \in W_{2}$ has pure weight $w \in \mathbb{Z}$ is $P$ is a linear combination of elements $\sigma^{\alpha}.\partial^{\beta}$ such that $w = w(\alpha) - w(\beta)$.

\begin{lemma}\label{weight}
A element $P$ in $W_{2}$ has pure weight $w$ if and only if it satisfies
\begin{equation}
U_{0}.P - P.U_{0} = w.P \quad {\rm or \ equivalently} \quad  U_{0}.P = P.(U_{0} + w)
\end{equation}
where $U_{0} := \sum_{h=1}^{k} \ h.\sigma_{h}.\partial_{h} .$
\end{lemma}

\parag{Proof} We leave to the reader to check the identity
$$ U_{0}.\sigma^{\alpha}.\partial^{\beta} - \sigma^{\alpha}.\partial^{\beta}.U_{0} = (w(\alpha) - w(\beta)).\sigma^{\alpha}.\partial^{\beta}.$$
Then if $P$ has pure weight $w$ it satisfies $(1)$.\\
Conversely, assume that $P \in W_{2}$ satisfies $(1)$. Then we may write in an unique way
$$ P = \sum_{v \in \mathbb{Z}}  \ P_{v} $$
where $P_{v}$ has pure weight $v$ and the sum is finite. This implies
$$ U_{0}.P = \sum_{v \in \mathbb{Z}} P_{v}.(U_{0} + v) = \sum_{v \in \mathbb{Z}} P_{v}.(U_{0} + w) .$$
So this implies that $\sum_{v \in \mathbb{Z}} (w - v).P_{v} = 0 $. By uniqueness of the decomposition in pure weight elements, we conclude that $P_{v} = 0$ for each $v \not= w$, and so $P = P_{w}$ concluding the proof.$\hfill \blacksquare$\\

\parag{Example} Define $U_{-1} := \sum_{h=0}^{k-1} (k-h).\sigma_{h}.\partial_{h+1}$. This is the element $T_q(S_{1})$ where $S_{1} \in W_{1}^{\mathfrak{S}_{k}}$ is the first symmetric function of $\partial_{x_{1}}, \dots, \partial_{x_{k}}$, so $S_{1} := \sum_{j=1}^{k} \partial_{x_{j}}$ and $q : \C^k \to \C^k$ the quotient map by the action of the action of $\mathfrak{S}_k$.\\
Then we have $U_{0}.U_{-1} - U_{-1}.U_{0} = - U_{-1}$ which corresponds to the fact that $U_{-1}$ has pure weight equal to $-1$.

\begin{lemma}\label{commutators}
We have the following commutation relations
\begin{align*}
& [U_{0}, A_{p, q}] = -(p+q).A_{p, q} \quad {\rm for} \ p, q, p+1, q-1 \in [1, k] \\
& [U_{0}, T^{h}] = -h.T^{h} \quad \forall h \in [2, k] \\
& [U_{0}, \tilde{T}^{h}] = -h.\tilde{T}^{h} \quad \forall h \in [2, k] \\
& [U_{-1}, A_{p, q}] = (k-p+1).A_{p+1, q-1} + (k-q).A_{p, q+1} \quad {\rm for} \ p, q, p+1, q-1 \in [1, k] \\
& [U_{-1}, T^{h}] = -(k-h).T^{h+1} - (k-1).A_{1, h} \quad \forall h \in [2, k] \\
&  [U_{-1}, \tilde{T}^{h}] = -(k-h).\tilde{T}^{h+1} - (k-1).A_{1, h} \quad \forall h \in [2, k]
\end{align*}
\end{lemma}

\parag{Proof} The first three commutation relations are obvious thanks to the previous lemma as the differential operators  $A_{p, q}$, $T^{h}$ and $\tilde{T}^{h}$ have pure weights respectively equal to $-(p+q), -h$ and $-h$.\\
This is not the case for the last two commutation relations which are given by the following elementary computations :
\begin{align*}
& [U_{-1}, \partial_{h}] = -(k-h).\partial_{h+1} \quad \forall h \in [1, k]   \quad {\rm for} \  h \in [2, k]  \ {\rm and \ we \ have} : \\
& [U_{-1}, \partial_{1}\partial_{h-1}] = U_{-1}\partial_{1}\partial_{h-1} - \partial_{1}\partial_{h-1}U_{-1}  = (\partial_{1}U_{-1} - (k-1).\partial_{2})\partial_{h-1} - \partial_{1}\partial_{h-1}U_{-1} \\
&   [U_{-1}, \partial_{1}\partial_{h-1}]  = \partial_{1}(\partial_{h-1}U_{-1} - (k-h+1).\partial_{h}) - (k-1).\partial_{2}\partial_{h-1} - \partial_{1}\partial_{h-1}U_{-1} \\
&  [U_{-1}, \partial_{1}\partial_{h-1}]   = -(k-h+1).\partial_{1}\partial_{h} - (k-1).\partial_{2}\partial_{h-1}
\end{align*}
Also, if $E := \sum_{h=1}^{k} \sigma_{h}.\partial_{h}$ we have
\begin{align*}
& [U_{-1}, \sigma_{h}.\partial_{h}] =  \sum_{p=0}^{k-1} (k-p).[\sigma_{p}.\partial_{p+1}, \sigma_{h}.\partial_{h}] \quad {\rm where} \ \sigma_{0} = 1 \\
& [U_{-1}, \sigma_{h}.\partial_{h}]  =  (k-h).[\sigma_{h}.\partial_{h+1}, \sigma_{h}.\partial_{h}] + (k-h+1).[\sigma_{h-1}.\partial_{h}, \sigma_{h}.\partial_{h}]  \quad {\rm then}\\
& [U_{-1}, E] = -\sum_{h=1}^{k} (k-h).\sigma_{h}.\partial_{h+1} + \sum_{h=1}^{k} (k-h+1).\sigma_{h-1}.\partial_{h}\\
&  [U_{-1}, E] = k.\partial_{1} \quad {\rm and \ then} \\
& [U_{-1}, \partial_{h} E] = U_{-1}\partial_{h}E - \partial_{h}EU_{-1} = (\partial_{h}U_{-1} -(k-h).\partial_{h+1}).E - \partial_{h}EU_{-1} \\
& [U_{-1}, \partial_{h} E]  = \partial_{h}.(EU_{-1} +k.\partial_{1}) - (k-h).\partial_{h+1}.E - \partial_{h}EU_{-1} \\
&  [U_{-1}, \partial_{h} E] = k.\partial_{1}\partial_{h} - (k-h).\partial_{h+1}E
\end{align*}
So we obtain for $h \in [2, k]$ :
\begin{align*}
& [U_{-1}, T^{h}] =  -(k-h+1).\partial_{1}\partial_{h} - (k-1).\partial_{2}\partial_{h-1} + k.\partial_{1}\partial_{h} - (k-h).\partial_{h+1}E \\
&  [U_{-1}, T^{h}] = -(k-h).T^{h+1} + (k-1).(\partial_{1}\partial_{h} - \partial_{2}\partial_{h-1}).
\end{align*}
Also
\begin{align*}
& [U_{-1}, \tilde{T}^{h}] =  [U_{-1}, T^{h}] + [U_{-1}, \partial_{h}] =  -(k-h).T^{h+1} + (k-1).(\partial_{1}\partial_{h} - \partial_{2}\partial_{h-1}) - (k-h).\partial_{h+1} \\
& [U_{-1}, \tilde{T}^{h}] = -(k-h).\tilde{T}^{h+1} + (k-1).A_{1, h} 
\end{align*}
concluding the proof.$\hfill \blacksquare$\\

The following corollary is obvious.

\begin{cor}\label{fond.2}
We have the inclusions
$$ \mathcal{I}.U_{0} \subset \mathcal{I}, \quad  \mathcal{I}.U_{-1} \subset \mathcal{I}, \quad  \tilde{\mathcal{I}}.U_{0} \subset \tilde{\mathcal{I}} \quad  {\rm and} \quad  \tilde{\mathcal{I}}.U_{-1} \subset \tilde{\mathcal{I}}. .$$
This implies that the right multiplications by $U_{0}$ and $U_{-1}$ define two endomorphisms of left $D_{N}-$modules on $\mathcal{M} = D_{N}\big/\mathcal{I}$ and $\tilde{\mathcal{M}} = D_{N}\big/\tilde{\mathcal{I}}$ .$\hfill \blacksquare$
\end{cor}

\parag{Remark} The commutation relation $[U_0, U_{-1}] = - U_{-1}$, which is consequence of the fact that $U_{-1}$ has pure weight $-1$, implies that we have a natural action on the $D_N-$modules $\mathcal{M}$ and $\tilde{\mathcal{M}}$ of the ring $\C< u, v >$ where the variables satisfies $(u+1).v = v.u$, where $u$ acts as $U_0$ and $v$ as $U_{-1}$. This $\C< u, v >-$structure then exists for instance on any sheaf ${Ext^q}_{D_N}(\mathcal{M}, \mathcal{N})$ or ${Ext^q}_{D_N}(\tilde{\mathcal{M}}, \mathcal{N})$  for any left  $D_N-$module $\mathcal{N}$ and any $q \in \mathbb{N}$. \\

\subsection{The theorem}

Recall that the derived Newton function introduced in \cite{[B.19]},  $DN_{m} := \tilde{T}(z^{m+k-1})$ for $m\geq -k+1$, is a polynomial of pure weight $m$ on $N$ (see also formula $(2)$ and $(2 bis)$ below). So they vanishes for $m$ in the interval  $[-k+1, -1]$.\\

 We shall use the following lemma later on.

\begin{lemma}\label{derive}
For each $m \geq 1$ we have the formula
$$U_{-1}[DN_{m}] = ( m+k-1).DN_{m-1}   .$$
\end{lemma}

\parag{Proof} We have, for $R \gg \vert \sigma\vert$
\begin{equation}
 DN_{m}(\sigma) = \frac{1}{2i\pi}\int_{\vert \zeta\vert = R} \frac{\zeta^{m+k-1}.d\zeta}{P_{\sigma}(\zeta)} .
 \end{equation}
Then
$$ \partial_{h}DN_{m} = -\frac{1}{2i\pi}\int_{\vert \zeta\vert = R} \ (-1)^{h}.\zeta^{k-h} \frac{\zeta^{m+k-1}.d\zeta}{P_{\sigma}(\zeta)^{2}} $$
and so
\begin{align*}
& U_{-1}[DN_{m}] =  \frac{1}{2i\pi}\int_{\vert \zeta\vert = R} \frac{\zeta^{m+k-1}.d\zeta}{P_{\sigma}(\zeta)^{2}}.P'_{\sigma}(\zeta) \\
&  U_{-1}[DN_{m}] =   -\frac{1}{2i\pi}\int_{\vert \zeta\vert = R} d\big(\frac{\zeta^{m+k-1}}{P_{\sigma}(\zeta)}\big) + (m+k-1). \frac{1}{2i\pi}\int_{\vert \zeta\vert = R} \frac{\zeta^{m+k-2}.d\zeta}{P_{\sigma}(\zeta)}\\
&  U_{-1}[DN_{m}] = (m+k-1).DN_{m-1} \qquad  \qquad \qquad \qquad \qquad \qquad \qquad \qquad \qquad \qquad \qquad \hfill \blacksquare\\
\end{align*}

We will show now the following characterization of global traces forms which completes the proof of the theorem \ref{solutions}.

\begin{thm}\label{caract.}
Let $G$ be a holomorphic function on $N$ which solution of the system $S\, 2)$ (or a global section of $Hom_{D_{N}}(\tilde{\mathcal{M}}, \mathcal{O}_{N})$). Then there exists a unique holomorphic function $g$ on $\C$ vanishing at order $k-1$ at the origin and such that $G = \tilde{T}(g)$.
\end{thm}

Recall that the equality $G = \tilde{T}(g)$ means that for each $\sigma$ such that $\Delta(\sigma) \not= 0$ we have 
\begin{equation*}
 G(\sigma) =  \tilde{T}(g)(\sigma) = \sum_{j=1}^{k} \frac{g(z_{j})}{P'_{\sigma}(z_{j})} \tag{2 bis}
 \end{equation*}
where $z_{1}, \dots, z_{k}$ are the roots of the polynomial $P_{\sigma}[z] = z^{k} + \sum_{h=1}^{k} (-1)^{h}.\sigma_{h}.z^{k-h}$ with discriminant $\Delta(\sigma)$.\\

We shall first prove the uniqueness statement of the theorem by computing the kernel of the (linear) map $\tilde{T} : \mathcal{O}(\C) \to \mathcal{O}(N)$.

\begin{lemma}\label{kernel}
The kernel of $\tilde{T}$ is the vector space of polynomials in $\C[z]$  of degree at most $k-2$.
\end{lemma}

\parag{Proof} Remark first that the vanishing of the derived Newton functions $DN_{m}$ for $m \in [-k+1, -1]$ implies that any polynomial of degree at most $k-2$ in $\C[z]$ is in the kernel of $\tilde{T}$.\\
Let $g \in \mathcal{O}(\C)$ be in the kernel of $\tilde{T}$. Write $g(z) := \sum_{0}^{+\infty} \lambda_{j}.z^{j+k-1}$. Then we have  $\tilde{T}(g) = \sum_{j=0}^{+\infty} \lambda_{j}.DN_{j}$ and this is the pure weight decomposition 
of $0$ so we have $\lambda_{j}.DN_{j} = 0$ for each $j \geq 0$. The only point to prove is the fact that for each $j \geq 0$ the polynomial $DN_{j}$ is not identically zero.  But we know that $N_{m}$ is a monic polynomial in $\sigma_{1}$ of degree $m$ and that $\partial_{1}N_{m}\big/m = DN_{m-1}$ for each $m \geq 0$ (see lemma \ref{derive}). So $DN_{j}$ is a monic polynomial in $\sigma_{1}$ of degree $j$. Then each $\lambda_{j}$ vanishes, and $g = 0$.$\hfill \blacksquare$\\

The proof of the theorem will use also the following lemma.

\begin{lemma}\label{poids}
For each integer $m \geq 1$ the kernel of $U_{-1}^{m}$ does not contain a element in $\C[\sigma_{1}, \dots, \sigma_{k}]$ of pure weight $m$.
\end{lemma}

\parag{Proof} For $m = 1$ a polynomial of pure weight $1$ is of the form $\lambda.\sigma_{1}$ with $\lambda \in \C$. Write $U_{-1} = k.\partial_{1} + V$ where the vector field $V$ satisfies $V[\sigma_{1}] = 0$. Then  we get $U_{-1}[\lambda.\sigma_{1}] = k.\lambda$ and this vanishes if and only if $\lambda = 0$ proving our assertion for $m = 1$.\\
Assume now that the lemma is proved for the integer $m-1 \geq 1$. We shall prove it for $m$.\\
Let $g$ a pure weight $m$ polynomial such that $U_{-1}^{m}[g] = 0$. Then, as $U_{-1}[g]$ has pure weight $m-1$ and is in the kernel of $U_{-1}^{m-1}$, the induction hypothesis gives $U_{-1}[g] = 0$.  Write $g = a.\sigma_{1}^{m} + h$ where $h$ has degree at most equal to $m-2$ in the variable $\sigma_{1}$ and where $a$ is a complex number. Now $U_{-1}[g] = 0$ implies 
$$ k.a.m.\sigma^{m-1} + k.\partial_{1}h + V[h] = 0.$$
As the degrees in $\sigma_{1}$ of $\partial_{1}h$ and of $V[h]$ are at most equal to $m-2$, we conclude that $a = 0$. This means that there exists an integer $p \geq 2$ such that $g$ may be written as
$g = u.\sigma_{1}^{m-p} + v $ where $u$ is in $\C[\sigma_{2}, \dots, \sigma_{k}]$ and has pure weight $p$ and where the degree of $v$ in $\sigma_{1}$ is at most $m-p-1$. Then the equality $U_{-1}[g] = 0$ implies
$$ k.(m-p).u.\sigma_{1}^{m-p-1} + k.\partial_{1}.v + V[u].\sigma_{1}^{m-p} + V[v] = 0 .$$
As the only term of degree $\geq m-p$ in $\sigma_{1}$ in the left hand-side is $V[u].\sigma_{1}^{m-p}$ we obtain  $V[u] = 0$;  and as $u$ does not depend on $\sigma_{1}$ we have $U_{-1}[u] = 0$. So $U_{-1}^{p}[u] = 0$ and as $u$ has pure weight $p \in [2, m]$, the induction hypothesis implies $u = 0$ if we have $p \not= m$. So the only case where we do not conclude that $g = 0$ is the case where $g$ does not depend on $\sigma_{1}$. For $k = 1$ $g$ is constant. For $k \geq 2$ the vanishing of $U_{-1}[g]$ implies that
 $$(k-1).\sigma_{1}.\partial_{2}g = - \sum_{h=2}^{k-1} \sigma_{h}.\partial_{h+1}g $$
 and the right hand-side is independent of $\sigma_{1}$ so $\partial_{2}g = 0$ and  $g$ does not depend on $\sigma_{1}$ and $\sigma_{2}$. So for $k = 2$ $g$ is constant. If $k \geq 3$  then 
 $$(k-2).\sigma_{2}.\partial_{3}g = - \sum_{h=3}^{k-1} (k-h).\sigma_{h}.\partial_{h+1}g$$
  and the right hand-side is independent on $\sigma_{2}$ so $\partial_{3}g = 0$ and {\it etc} ... \\
 We conclude that $g$ is constant and of pure weight $m \geq 2$, ending the proof of the induction step.$\hfill \blacksquare$\\

\parag{Proof of the theorem \ref{caract.}} Consider a holomorphic function $G$ on $N$ which is solution of the system $S\, 2)$ and write
$$ G = \sum_{w \geq 0} G_{w} $$
where $G_{w}$ is a pure weight $w$ polynomial. This series converges uniformly on any compact set in $N$ and as the differential operators which generate $S\, 2)$ have pure weights, for each $w$ the polynomial $G_{w}$ is a solution of $S\, 2)$. So it is enough to prove the theorem when we assume that $G$ is a polynomial of pure weight $w$.\\
Assuming now that $G = G_{w}$ we have $U_{-1}^{w}[G_w]$ which is still a solution of $S\, 2)$, thanks to the corollary \ref{fond.2},  and has pure weight $0$. So it is a constant $\lambda \in \C$. Consider now the polynomial $G_{w} - \mu.DN_{w}$ where $\mu \in \C$ is chosen in order that $U_{-1}^{w}(G_{w} - \mu.DN_{w})= 0$. This is possible thanks to the lemma \ref{derive}. and the fact that $DN_w$ is a trace form. Then the lemma \ref{poids} implies that $G_{w} = \mu.DN_{w}$.$\hfill \blacksquare$\\

\parag{Remark} The theorem \ref{caract.} extends immediately to a holomorphic function on a connected open set containing the origin (with the same proof), because we can again use the Taylor expansion at the origin in order to reduced the question to pure weight polynomials as above; and then to conclude by analytic continuation.\\

\subsection{Proofs of the corollaries}

\parag{Proof of  corollary \ref{charact. bis}} We shall show first that $\partial_{1}F$ is a trace form. Thanks to theorem \ref{caract.} it is enough to prove that $G := \partial_{1}F$ satisfies the system $S\, 2)$.\\
Remark that for $p, q$ with $p \in [1, k-1]$ and $q \in [2, k]$ we have
$$ \partial_{q}T^{p+1} - \partial_{p+1}.T^{q} = \partial_{1}.(\partial_{p}\partial_{q} -\partial_{p+1}\partial_{q-1}) = A_{p, q}.\partial_{1} .$$
This implies that $A_{p, q}[G] = 0$ for each $p, q$ such that $p, p+1, q, q-1$ belongs to $[1, k]$.\\
The commutation relations $[E, \partial_{h}] = - \partial_{h}$ implies 
$$ \tilde{T}^{h} \partial_{1} = \partial_{1}T^{h} \quad \forall h \in [2, k]$$
showing that $\tilde{T}^{h}[G] = 0$ for each $h \in [2, k]$. So $G$ satisfies $S\, 2)$ and then is a trace form. Let $g \in \mathcal{O}(\C)$ satisfies for each $\sigma \in N$
$$ G(\sigma) = \tilde{T}(\gamma)(\sigma) \quad {\rm where} \quad \gamma(z) =  z^{k-1}.g(z) ,$$
and let $F_{0} := T(f)$ with $f' := g$.

Then we obtain, using the lemma below, that  $\partial_{1}(F_{0} -F) = 0$ . Now $T^{h}[F_{0}- F] = 0$ for each $h \in [2, k]$ implies, as $\partial_{1}(F_{0}-F) = 0$ and $\mathcal{T}^h = \partial_1\partial_{h-1} + \partial_hE$, that $\partial_{h}E(F_{0}-F) = 0$  for each $h \in [2, k]$. Then $E(F_{0}-F) $ is constant because $F_{0}- F$ is independent of $\sigma_{1}$, and then vanishes.  So $ F = F_{0}+ c$ where $c$ is a constant. $\hfill \blacksquare$\\

\begin{lemma}  For each $f \in \mathcal{O}(\C)$ we have the identity
$$\partial_{1}T(f) = \tilde{T}(z^{k-1}.f'(z)).$$
\end{lemma}

\parag{Proof} We have, for $R \gg 1$ (see formula $(M)$ or \cite{[B.19]})
$$ T(f)(\sigma) = \frac{-1}{2i\pi}\int_{\vert \zeta\vert = R} f'(\zeta).Log\big(P_{\sigma}(\zeta)\big/\zeta^{k}\big).d\zeta + k.f(0) $$
so $$ \partial_{1}T(f)(\sigma) = \frac{1}{2i\pi}\int_{\vert \zeta\vert = R} f'(\zeta)\frac{\zeta^{k-1}}{P_{\sigma}(\zeta)}.d\zeta = \tilde{T}(z^{k-1}.f'(z))$$
concluding the proof.$\hfill \blacksquare$\\

 \parag{Proof of corollary  \ref{3 to 1}} Using the point $8.$ in the proposition \ref{direct relations}, let $F$ be the entire holomorphic function such that $\partial_hF = (-1)^{k-h-1}\Phi_{k-h}$ for $h \in [1, k]$. As the components of $\Phi$ are solutions of $S\, 2)$ the relations of commutation show that $T^mF$ and $A_{p, q}F$ are constant. Let $F_w$ the pure weight $w$ part of $F$. Recall that the Taylor series converges uniformly on each compact set in $N$ and this gives a pure weight decomposition $F = \sum_{w = 0}^{\infty} F_w$. As $T^m$ and $\tilde{T}^m$ have pure weight $-m$, for each  $m \in [2, k]$, the uniqueness of this decomposition implies  $\partial_hT^mF_w = 0, \forall h \in [1,k]$. This implies that $T^mF_w$  are constant for each $m \in [2, k]$ and $w \in \mathbb{N}$. But $T^mF_w$ has pure weight $w - m$ and a non zero constant has pure weight $0$. So for $w \geq k+1$ we conclude that $T^mF_w = 0$ for each $m \in [2, k]$. Then the corollary \ref{charact. bis} implies that $F_w$ is a trace function for $w \geq k+1$. For $w \in [0, k]$ we have $T^mF_w = 0$ for $w \not= m$ and so we see that there exists a constant $c_w$ such that $T^m(F_w -c_w.\sigma_w ) = 0$ for all $m \in [2, k]$. We conclude that $F - \sum_{w = 1}^k c_w.\sigma_w$ is a trace function, thanks to the corollary \ref{charact. bis}. \\
 As $\sigma_1 = N_1$ is a trace function, $G := F -  \sum_{w = 2}^k c_w.\sigma_w$ is a trace function, and using the theorem \ref{solutions} and the point 1. of the proposition \ref{direct relations} we see that the constant vector $C$ with components $(-1)^h.c_{k-h}$ for $h \in [2, k]$ and $0$ for $h = k-1$ is a solution of $S\, 3)$.\\
 To complete the proof we shall use the following lemma.

\begin{lemma}\label{cste sol.} 
Let $\Phi$ be a constant vector in $\C^k$ which is solution of $S\, 3)$. Then $\Phi$ is equal to $c.V$ where $V$ is the vector $V := \begin{pmatrix} 0 \\ .\\ . \\ 0 \\ 1\end{pmatrix}$. 
\end{lemma}

\parag{Proof} So we have the relations $\partial_k(A^p\Phi) = 0$ for each $p \in [1, k-1]$. As $^tA$, the transpose of $A$,  is the matrix of the multiplication by $z$ in the $\C[\sigma]-$basis $1, z, \dots, z^{k-1}$  of the $\C[\sigma]$ algebra $R := \C[\sigma, z]\big/(P_\sigma(z))$, defining the derivation $\partial_k$  on $R$ by $\partial_k(z) = 0, \partial_k(\sigma_h) = 0$, for $h \not= k$ and $\partial_k(\sigma_k) = 1$, our problem is equivalent to prove the following\footnote{warning : the transposition reverse here the numbering of the components }:
\begin{itemize}
\item Let $Q \in \C[z] $ be a polynomial of degree $q$ at most equal to $k-1$ and let $[Q]$ be the element in $R$ define by $Q$. Assume that $\partial_k(z^p[Q](z)) = 0$ for each $p \in [1, k-1]$. Then $[Q]$ is a constant in $\C \subset R$.
\end{itemize}
We shall prove this assertion by induction on the degree $q$ of $Q$. First remark that if $Q$ is the constant $c$ we have $\partial_k(z^p.[c]) = 0$ for each $p \in [1, k-1]$.\\
Assume that for $Q$ of degree $q-1$, with $1 \leq q \leq k-1$ the assertion is proved. If $Q$ has degree $q$ with dominant coefficient $\gamma.z^q$, consider the class of $z^{k-q}.Q(z)$ in $R$. The only  coefficient in the basis $1, z, \dots, z^{k-1}$ of $[z^{k-q}.Q(z)]$ where $\sigma_k$ appears  is the coefficient of $1$ which is equal to $(-1)^{k-1}.\gamma.\sigma_k$. So, as $k-q$ is in $[1, k-1]$, the condition $\partial_k([z^{k-q}.Q(z)]) = 0$ implies $\gamma = 0$ and we conclude that $Q$ is constant by the induction hypothesis. This conclude the proof.$\hfill \blacksquare$\\

\parag{End of the proof of \ref{3 to 1}} So there exists $g \in \mathcal{O}(\C)$ such that $F = T(g)$. And we have
$$\partial_hT(g) = (-1)^{h-1}.V\tilde{T}(g')_{k-h} = (-1)^{k-h-1}.\Phi_{k-h} = (-1)^{k-h-1}.V\tilde{T}((-1)^k.f)_{k-h}.$$
The linear map $f \mapsto V\tilde{T}(f)$ is injective\footnote{In term of trace forms this is consequence of the non vanishing of the Van der Mond determinant in $z_1, \dots, z_k$ when $\Delta(\sigma) \not= 0$. See also prop. 2.3 in \cite{[B-MF]}.} so we conclude that $g' = (-1)^k.f$.$\hfill \blacksquare$\\

\subsection{Lagrange interpolation and Lisbon integrals}

We begin by recalling the Lagrange interpolation formula.

\begin{lemma}\label{interpol.0}
Let $f $ be an entire holomorphic function on $\C$ and let $\sigma \in N$. Then the Lagrange interpolation polynomial of $f$ for the monic polynomial $P_\sigma$ is given by the following formula
\begin{equation}
\Pi_f(\sigma)[z] := \frac{1}{2i\pi}\int_{\vert \zeta \vert = R} f(\zeta).\frac{P_\sigma(\zeta) - P_\sigma(z)}{(\zeta - z).P_\sigma(\zeta)}.d\zeta 
\end{equation}
for $R > \vert z\vert$ and $R \gg \vert\vert \sigma \vert\vert$.
\end{lemma}

\parag{Proof}  We have
$$ \Pi_f(\sigma)[z] = -\frac{P_\sigma(z)}{2i\pi}\int_{\vert \zeta \vert = R} \frac{f(\zeta)}{(\zeta - z).P_\sigma(\zeta)}.d\zeta +  \frac{1}{2i\pi}\int_{\vert \zeta \vert = R} \frac{f(\zeta)}{\zeta - z}.d\zeta $$

and the residue formula gives
\begin{equation*}
\Pi_f(\sigma)[z] = f(z) - P_\sigma(z).Q_f(\sigma, z) \quad {\rm with} \quad  Q_f(\sigma, z) = \frac{1}{2i\pi}\int_{\vert \zeta \vert = R} \frac{f(\zeta)}{(\zeta - z).P_\sigma(\zeta)}d\zeta 
\end{equation*}
where $Q_f(\sigma, z)$ is holomorphic on $N \times \C$. The relation 
\begin{equation}
f(z) = \Pi_f(\sigma)[z] + P_\sigma(z).Q_f(\sigma, z)
\end{equation}
 is then clearly the euclidian division of $f$ as an entire function on $N \times \C$ by the monic polynomial $P_\sigma(z)$ of degree $k$, as the formula defining $\Pi_f$ shows that it is a degree at most $k-1$ polynomial in $z$ with holomorphic coefficients in $\mathcal{O}(N)$.\\
So $\Pi_f$ is the Lagrange interpolation polynomial of $f$ for $P_\sigma$.$\hfill \blacksquare$\\

\begin{cor}\label{interpol.1}
Let $\Phi$ be the vector Lisbon integral associated to  $f \in \mathcal{O}(\C)$. Then the coefficient of $z^h$ in  $\Pi_f(\sigma)$ is given by the following formula (where $\sigma_0 \equiv 1$)
\begin{equation}
\Pi_f(\sigma)_h =  \sum_{p = 0}^{k-h-1} (-1)^p.\sigma_p.\Phi_{k-p-h-1} \quad \forall h \in [0, k-1]
\end{equation}
\end{cor}

\parag{Proof} First let us compute the coefficient $\gamma_h(\zeta)$ of $z^h$ in the polynomial
 $$\frac{P_\sigma(\zeta) - P_\sigma(z)}{\zeta - z}$$
  in $\C[\zeta, z] $. As for $m \in \mathbb{N}^*$ we have
$$ \frac{\zeta^m - z^m}{\zeta - z} = \sum_{h=0}^{m-1} \  \zeta^{m-h-1}.z^h $$
we obtain that 
$$ \gamma_h(\zeta) = \sum_{p = 0}^{k-h-1} (-1)^p.\sigma_p.\zeta^{k-p-h-1} .$$
This implies the formula $(5)$.$\hfill \blacksquare$\\

\parag{Remark} Remark that  the linear system linking $\Phi$ and $\Pi_f$ as elements in $\mathcal{O}(N)^k$ is triangular with diagonal elements equal to $1$ and  the corresponding matrix does not involve $\sigma_k$.  In particular $\Phi_0$ is the coefficient of $z^{k-1}$ in $\Pi_f$ and the same linear system links the vectors $\partial_k(\Phi)$ and $\partial_k(\Pi_f)$.\\


\providecommand{\bysame}{\leavevmode\hbox to3em{\hrulefill}\thinspace}
\providecommand{\MR}{\relax\ifhmode\unskip\space\fi MR }
\providecommand{\MRhref}[2]{%
  \href{http://www.ams.org/mathscinet-getitem?mr=#1}{#2}
}
\providecommand{\href}[2]{#2}
\begin{thebibliography}{}

\end{thebibliography}


\bigskip

\begin{thebibliography}{99}
 
 \parag{References}

 
 \bibitem{[B.19]} Barlet, D. {\it On Symmetric Differential Operators}: math-arXiv: 1911.09347
 
 \bibitem{[B-MF]} Barlet, D. and Monteiro-Fernandes, T. {\it On Lisbon integrals}: math-arXiv: 1906.09801 to appear in Math. Z.
 
 \bibitem{[Bjork]} Bjork, J.E. {\it Analytic $D-$modules and Applications} Springer (1993).
 

 
 \end{thebibliography}
\end{document}